\documentclass[reqno,12pt]{amsart}
\usepackage{amssymb}
\usepackage{amsmath} %new
\usepackage{epsf} %new
\usepackage{latexsym} %new
\usepackage{amscd}
\usepackage{graphics}
\usepackage{graphpap}
\usepackage{soul}
\usepackage[all]{xypic}

\def\CC {{\mathbb C}}

\def\NN {{\mathbb N}}

\def\QQ {{\mathbb Q}}
\def\RR {{\mathbb R}}
\def\TT {{\mathbb T}}
\def\XX {{\mathbb X}}

\def\ZZ {{\mathbb Z}}
\def\sA {{\mathcal A}}
\def\sB {{\mathcal B}}
\def\sC {{\mathcal C}}

\def\sK {{\mathcal K}}
\def\sL {{\mathcal L}}

\def\sQ {{\mathcal Q}}

\def\sS {{\mathcal S}}

\def\sV {{\mathcal V}}

\def\Hom {\mathrm{Hom}}

\def\gl {\lambda}
\def\supp {\mathrm{supp}}

\def\and {\wedge}
\def\qed {{\blacksquare}}
\def\onto {\,\rule[.04in]{.15in}{.01in}\kern-4pt\raise.3ex\hbox{$\scriptscriptstyle>\!>\,$}}
\def\to {\longrightarrow}

\def\st {such that}

\def\sii {if and only if}

\def\< {{\langle}}
\def\> {{\rangle}}

\def\gb {{\beta}}

\def\gl {{\lambda}}

\def\cg {{\overline{g}}}
\def\cx {{\overline{x}}}
\def\ct {{\overline{t}}}
\def\^ {{\widehat{ \:}}}

\def\pf {{\em \noindent Demostraci\'on:}}
\def\epf {~\hfill $\qed$}

\def\pf {{\em \noindent Proof:}} %new
\def\l( {{\left)}} %new
\def\r( {{\right)}} %new
\def\l[ {{\left[}} %new
\def\r] {{\right]}} %new
\def\l{ {{\left{}} %new
\def\r} {{\right}}} %new

\def\ha{\widehat{\phantom{x}}}
\def\haha{\hat{\phantom{t}}\hat{\phantom{t}}}

\newcommand{\mkp}{\medskip}

\newfont{\cyr}{wncyr8}
\newfont{\cyb}{wncyr8}
\usepackage{pdfpages}
\usepackage{amsxtra}
\newtheorem{thm}{Theorem}[section]
\newcommand{\bthm}{\begin{thm}}
\newcommand{\ethm}{\end{thm}}

\newtheorem{prop}[thm]{Proposition}
\newcommand{\bprp}{\begin{prop}}
\newcommand{\eprp}{\end{prop}}

\newtheorem{fact}[thm]{Fact}
\newcommand{\bfct}{\begin{fact}}
\newcommand{\efct}{\end{fact}}

\newtheorem{prob}[thm]{Problem}
\newcommand{\bprb}{\begin{prob}}
\newcommand{\eprb}{\end{prob}}

\newtheorem{quest}[thm]{Question}
\newcommand{\bqtn}{\begin{quest}}
\newcommand{\eqtn}{\end{quest}}

\newtheorem{lem}[thm]{Lemma}
\newcommand{\blem}{\begin{lem}}
\newcommand{\elem}{\end{lem}}

\newtheorem{claim}[thm]{Claim}
\newcommand{\bclm}{\begin{claim}}
\newcommand{\eclm}{\end{claim}}

\newtheorem{cor}[thm]{Corollary}
\newcommand{\bcor}{\begin{cor}}
\newcommand{\ecor}{\end{cor}}

\newtheorem{conj}[thm]{Conjecture}
\newcommand{\bcnj}{\begin{conj}}
\newcommand{\ecnj}{\end{conj}}

\theoremstyle{definition}
\newtheorem{defn}[thm]{Definition}
\newcommand{\bdfn}{\begin{defn}}
\newcommand{\edfn}{\end{defn}}

\newtheorem{spec}[thm]{Specializing}
\newcommand{\bspc}{\begin{spec}}
\newcommand{\espc}{\end{spec}}

\theoremstyle{remark}
\newtheorem{rem}[thm]{Remark}
\newcommand{\brem}{\begin{rem}}
\newcommand{\erem}{\end{rem}}

\newtheorem{cnv}[thm]{Convention}
\newcommand{\bcnv}{\begin{cnv}}
\newcommand{\ecnv}{\end{cnv}}

\newtheorem{exam}[thm]{Example}
\newcommand{\bexm}{\begin{exam}}
\newcommand{\eexm}{\end{exam}}

\newtheorem{exercise}[thm]{Exercise}
\newcommand{\bexr}{\begin{exercise}}
\newcommand{\eexr}{\end{exercise}}

\newtheorem{thmy}{\textbf{Theorem}}
\newenvironment{thmx}{\stepcounter{thm}\begin{thmy}}{\end{thmy}}

\newcommand{\X}{\mathbb X}

\renewcommand{\r }{\rangle}
\renewcommand{\l }{\langle}

\begin{document}
\label{begin-art}

\title{Tensor products of topological abelian groups and Pontryagin duality}

\date{01/27/2024}

\author[M. Ferrer]{Mar\'ia V. Ferrer}
\address{Universitat Jaume I, Instituto de Matem\'aticas de Castell\'on,
Campus de Riu Sec, 12071 Castell\'{o}n, Spain.}
\email{mferrer@mat.uji.es}

\author[J. Hern\'andez-Arzusa]{Julio Hern\'andez-Arzusa}
\address{Universidad de Cartagena, Departamento de Matem\'{a}ticas,
Campus de San Pablo, Cartagena, Colombia.}
\email{jhernandeza2@unicartagena.edu.co}

\author[S. Hern\'andez]{Salvador Hern\'andez}
\address{Universitat Jaume I, Departamento de Matem\'{a}ticas,
Campus de Riu Sec, 12071 Castell\'{o}n, Spain.}
\email{hernande@mat.uji.es}

\thanks{ The first and third authors acknowledge partial support by by the Spanish Ministerio de Econom\'{i}a y Competitividad,
grant: MTM/PID2019-106529GB-I00 (AEI/FEDER, EU) and by the Universitat Jaume I, grant UJI-B2022-39}

\begin{abstract}
Let $G$ be the group of all $\ZZ$-valued homomorphisms of the Baer-Specker group $\ZZ^\NN$. The group $G$ is algebraically
isomorphic to $\ZZ^{(\NN)}$, the infinite direct sum of the group of integers, and equipped with the topology of pointwise convergence on $\ZZ^\NN$, becomes
a non reflexive prodiscrete group. It was an open question to find its dual group $\widehat{G}$. Here, we answer this question by proving that
$\widehat{G}$ is topologically isomorphic to $\ZZ^\NN\otimes_\mathcal{Q}\TT$, the (locally quasi-convex) tensor product of $\ZZ^\NN$ and $\TT$.
Furthermore, we investigate the reflexivity properties of the { groups $C_p(X,\ZZ)$},
the group of all $\ZZ$-valued continuous functions on $X$ equipped with the pointwise convergence topology,
and $A_p(X)$, the free abelian group on a $0$-dimensional space $X$ equipped with
the topology $t_p(C(X,\ZZ))$ of pointwise convergence topology on $C(X,\ZZ)$. In particular, we prove that
$\widehat{A_p(X)}\simeq C_p(X,\ZZ)\otimes_\mathcal{Q}\TT$ and we establish the existence of $0$-dimensional spaces $X$ such that
$C_p(X,\ZZ)$ is Pontryagin reflexive.
\end{abstract}


%\dedicatory{}

\maketitle \setlength{\baselineskip}{24pt}

\section{Introduction}
In this paper, we are concerned with the duality properties of some \emph{prodiscrete} abelian groups, namely, projective limits of abelian discrete groups or,
equivalently, closed subgroups of products of abelian discrete groups, and \emph{proto-discrete} groups, that is, subgroups of products of abelian discrete groups
(cf. { \cite[Chapter 4, pag. 153]{HM-Pro-Lie:2023}, \cite[Cor. 3.13]{Hofmann_Morris:MZ2004})}.
{ Specifically}, we investigate the reflexivity properties of $C_p(X,\ZZ)$,
the group of all $\ZZ$-valued continuous functions on a $0$-dimensional topological space $X$, equipped with the pointwise convergence to\-po\-lo\-gy,
and $A_p(X)$, the free abelian group on a $0$-dimensional space $X$, equipped with
the to\-po\-lo\-gy $t_p(C(X,\ZZ))$ of pointwise convergence on $C(X,\ZZ)$. Our approach is based on the notion of \emph{locally quasi-convex tensor product} of
topological abelian groups. Using this notion we have been able { to find} the dual groups of $C_p(X,\ZZ)$ and $A_p(X)$ for a $0$-dimensional Hausdorff
topological space $X$.

All spaces are assumed to be Tikhonov and all groups are assumed to be abelian. If $G$ is a topological group,
a {\it character\/} of $G$ is a continuous homomorphism of $G$ to the
group $\Bbb{T}=\{z\in\CC:|z|=1\}$. For every topological group
$G$ the {\it dual group\/} $\widehat{G}$ is the topological group of all characters of $G$, equipped with
the compact-open topology. There is a natural evaluation homomorphism (not necessarily continuous)
$$\mathcal E :G\longrightarrow \widehat{\widehat G\,}$$ of $G$ to its bidual group, defined by $\mathcal E(g)(\chi):=\chi(g)$ for all $g\in G,\chi\in\widehat{G}$.
A topological group $G$ is {\it Pontryagin reflexive\/}, or {\it $P$-reflexive\/} for short,
if the evaluation map  $\mathcal E$ is a topological isomorphism between
$G$ and $\widehat{\widehat G\,}$. If $\mathcal E$ is just bijective
(but not necessarily continuous), we say that $G$ is {\it $P$-semireflexive}.

The Pontryagin--van Kampen theorem, saying that every locally compact abelian
group is $P$-reflexive, has been generalized to other classes of topological
groups. Kaplan \cite{kaplan1,kaplan2}
proved that the class of $P$-reflexive groups is closed
under arbitrary products and direct sums,
Smith \cite{smith} proved that Banach spaces
and reflexive locally convex vector spaces are $P$-reflexive
topological groups. Extending further the scope of the Pontryagin-van Kampen theorem, Banaszczyk \cite{Banaszczyk} introduced
the class of nuclear groups and investigated its duality properties.
The class of nuclear groups contains all locally compact abelian groups
and is closed under taking products, subgroups and quotients \cite[7.5 and 7.10]{Banaszczyk}. It follows from this that every proto-discrete abelian group is nuclear.

Hofmann and Morris have studied in depth the topological and algebraic structure of { pro-Lie} groups, that is, projective limits of Lie groups.
One essential ingredient in this setting are prodiscrete groups. In the Abelian case, a main goal is to develop a satisfactory duality theory
of prodiscrete groups (see \cite{HM-TP,HM-Axioms}).
Here, we further develop the research initiated in \cite{FerHerSha:pams} where the duality properties of different classes of prodiscrete groups are studied.
Thereby, we answer some questions that had been left open in the aforegoing paper.
We now formulate our main results:

\begin{thmx}\label{Thm_A}
Let $G$ denote the group $\ZZ^{(\NN)}\simeq\Hom(\ZZ^\NN,\ZZ)$, equipped with the topology $t_p(\ZZ^\NN)$ (inherited from $\ZZ^{\ZZ^\NN}$). Then we have
\emph{$\widehat{G}\cong \ZZ^\NN\otimes_\mathcal{Q} \TT$.}
\end{thmx}
\mkp

We notice that this result answers a question asked by Kazuhiro Kawamura
at the Conference on Set-Theoretic Topology and its Applications that took place in Yokohama (Japan) in 2015.

Our second main result extends Theorem A to a wider setting.

\begin{thmx}\label{Thm_B}
Let $X$ be a $0$-dimensional, $k_\NN$-space. Then we have
$$\widehat{A_p(X)}\simeq C_p(X,\ZZ)\otimes_\mathcal{Q}\TT.$$
\end{thmx}\mkp

The group $C_p(X,\ZZ)$ is not Pontryagin reflexive in general. For example, if $X$ is countable and nondiscrete then
$C_p(X,\ZZ)$ is metrizable but not complete. Thus, by a well known result of Au\ss enhofer \cite{aus} and Chasco \cite{chasco},
it follows that $C_p(X,\ZZ)$ is not Pontryagin reflexive. Nevertheless, there are $0$-dimensional spaces
$X$ such that $C_p(X,\ZZ)$ is Pontryagin reflexive. Indeed, in \cite{Her_Usp:JMAA2000}, it was proved the existence
of $0$-dimensional spaces for which the group $C_p(X,\RR)$, of all real-valued continuous functions on $X$, equipped with the
pointwise convergence topology, is Pontryagin reflexive. Furthermore, it was shown there that the answer to the question whether
$C_p(X,\RR)$ is Pontryagin reflexive depends on the axioms of set theory that we assume.
Going on further in this research line, we get:

\begin{thmx}\label{Thm_C}
{ There are nondiscrete $0$-dimensional spaces $X$} such that the groups $C_p(X,\ZZ)$ are
Pontryagin reflexive.
\end{thmx}

\section{Tensor products in the category of topological abelian groups}

In this section, we collect some basic facts about tensor products and reflective subcategories of topological abelian groups.
The tensor product of locally compact abelian (LCA) groups was introduced and investigated by Hofmann in \cite{Hofmann:1964},
where he approached it for any concrete subcategory $\sB$ of the category $\sL$ of locally compact abelian groups.
His definition is based on a \emph{natural} universal mapping principle but there is no general construction for the tensor product as in the { purely algebraic} case.
In this setting, Hofmann proved the existence and  unicity of the tensor product in some particular cases. For instance,
if we take for $\sB$ the subcategory $\sK$ of all compact abelian groups, then the $\sK$-tensor product of two compact abelian groups $G$ and $H$ exists
and is topologically isomorphic to the character group of the discrete group of all \emph{bimorphisms} (defined below) of $G\times H$ into the $1$-dimensional torus.
%However, although the $\sL$-tensor product of two locally compact abelian groups can be determined explicitly for some special cases,
%it cannot be established the existence of the tensor product for each specific subcategory $\sB$ of $\sL$ in general.

A different strategy was considered by Garling \cite{Garling}, who made use of the duality theory of locally compact abelian groups in order to
define the tensor product as a \emph{sort of dual group}. In this approach, one needs all the machinery of Pontryagin duality in order to define the
topological tensor product of two abelian groups. Although this approach works well for compactly generated locally compact abelian groups,
it is not so successful for general locally compact groups and it seems { quite unlikely} that this approach %to extend for tensor products of (finite) families
%of locally compact abelian groups. Furthermore, it might not
works so well for other categories of topological abelian groups not possessing a satisfactory duality theory.
For instance, the following question remains open for the category $\sA$ of topological abelian groups:
given a subcategory $\sB$ of $\sA$, is there a \emph{canonical} subcategory $\sB^\otimes\subseteq \sA$ such that
$\sB\subseteq \sB^\otimes$ and, for all $G,H\in\sB$, there is a tensor product $G\otimes_{\sB^\otimes} H\in\sB^\otimes$?
Garling solves this question for what he calls \emph{groups with duality} by proving that if $\mathcal{B}$ designates the subcategory of groups with duality
then $\sB = \sB^\otimes$ (incidently, this subcategory coincides with the subcategory $\sQ$ of (Hausdorff) \emph{locally quasi-convex groups} \cite{Vilenkin}).
In this setting, one may ask whether the subcategory $\sQ$ of locally quasi-convex groups is the smallest subcategory containing $\sL$ where the existence
of tensor product can be established.

Here, our approach follows that of Hofmann to investigate tensor products in different subcategories of topological abelian groups. As we have mentioned above, this methodology has the disadvantage that the tensor product does not always exist for each subcategory of topological abelian groups. However, we show that it does exist and it is unique,
up to topological isomorphisms, for large  subcategories of topological abelian groups.

%Here, we are concerned with these questions and our aim is twofold. First, we provide an affirmative answer to the question formulated above.
%As an application we obtain that $\sL^\otimes$ is the subcategory of projective limits of Lie abelian groups (the so called \emph{pro-Lie} abelian groups).
%Secondly, we apply our results on tensor products in order to identify the dual of the group $\bigoplus\limits_{n\in\ZZ} \ZZ$ equipped with a group topology that makes it prodiscrete (therefore pro-Lie), non discrete and non reflexive (see \cite{FerHerSha:pams}).

\bdfn\label{def1}
Let $\sA$ be the category of all topological Abelian groups and let $\sB$ be a subcategory of $\sA$. Given arbitrary $G,H\in\sB$ and $A\in\sA$,
the map $b\colon G\times H\to~A$ %($G,H,A\in\sA$)
is a \emph{continuous bihomomorphism}
if $x\rightarrow b(x,h)$ and $y\rightarrow b(g,y)$ are continuous group homomorphisms whenever $g\in G$ and
$h\in H$ are kept fixed, and if $b$ is continuous at $(0,0)$.
When $A$ is the $1$-dimensional torus $\TT$, we say that $b$ is a \emph{continuous bicharacter} of $G\times H$.

Given $G,H\in\sB$, the pair $(T, \otimes_{\mathcal{B}})$, with $T\in\mathcal{B}$  and $\otimes_{\mathcal{B}} \colon G\times H\longrightarrow  T$ a continuous  bihomomorphism,
is called a \emph{$\sB$-tensor product} of $G$ and $H$ if for every continuous { bihomomorphism} $b\colon G\times H\longrightarrow B$ whith $B\in \mathcal{B}$,
there exists a continuous homomorphism $\widetilde{b}\colon T\longrightarrow B$ such that $\widetilde{b}\circ \otimes_{\mathcal{B}}=b$. That is, the diagram
\vspace{-0.5cm}
\begin{equation*}
\begin{picture}(148,120)
\put(40,5){$B$}
\put(-55,60){$G\times H$}
\put(112,60){$T$}
\put(-15,60){\vector(1,0){121}}
\put(40,65){$\otimes_{\mathcal{B}}$}
\put(-21,56){\vector(4,-3){47}}
\put(-15,30){$b$}
\put(108,56){\vector(-4,-3){47}}
\put(100,30){$\widetilde{b}$}
\end{picture}
\end{equation*}
commutes. If there is a $\sB$-tensor product of any two groups $G,H$ in $\sB$, we say that \emph{$\sB$ has tensor product}.
If for a subcategory $\sB'\supseteq \sB$, we have that there is a $\sB'$-tensor product of any two groups $G,H$ in $\sB$,
we say that \emph{$\sB$ has tensor product in $\sB'$}.
\edfn
\mkp

\brem\label{Rem00}
From here on, the term \emph{group} means topological abelian group.
\erem

The following result is folklore. We include its proof here for completeness' sake.

\bthm\label{th 1}
The category $\sA$ of all topological abelian groups has tensor pro\-duct. Furthermore, the homomorphism $\widetilde{b}$ referred to in Definition \ref{def1} is unique.
\ethm
\pf\
  Let $G,H$ be two arbitrary groups in $\sA$. Then we first take $(A(G\times H),\sigma)$ be the free topological abelian group associated to $G\times H$,
where $\sigma\colon  G\times H \longrightarrow A(G\times H)$  is the canonical embedding.  Denote by $N$ %the closure in $A(G\times H)$ of
the subgroup generated { by} the elements of the form
\begin{center}
$\sigma(a+b,c)-\sigma(a,c)-\sigma(b,c)$\\$\sigma(a,c+d)-\sigma(a,c)-\sigma(a,d),$
\end{center}
where $a,b\in G$ and $c,d\in H$. Define $G\otimes H:=A(G\times H)/N$ and let $$\pi \colon A(G\times~H)\longrightarrow A(G\times H)/N$$ be
the canonical quotient map. We have that $\pi\circ\sigma$ is a continuous bihomomorphism by the definition of $N$.
The pair $(G\otimes H,\pi\circ\sigma)$ satisfies the requirements of being a tensor product.
Indeed, if  $b\colon G\times H\longrightarrow A$ is a continuous bihomomorphism, %where $A\in \mathcal{A}$, %then there exists a continuous homomorphism
%$h\colon A(G\times H)/N\longrightarrow A$ such that $h\circ\pi\circ\sigma =b$. Indeed,
%from the definition of $(A(G\times H),\sigma)$, there exists
then there exists a unique continuous homomorphism $\hat{b}\colon A(G\times H)\longrightarrow A$
such that $\hat{b}\circ \sigma=b$. Furthermore, since $N\subseteq \ker \hat{b}$, there exists a continuous homomorphism
$\widetilde{b}\colon A(G\times H)/N\longrightarrow A$ such that $\widetilde{b}\circ\pi =\hat{b}$.
This implies that $\widetilde{b}\circ (\pi\circ\sigma)=b$.
The uniqueness of $\tilde{b}$ follows from the fact that $\pi\circ\sigma(G\times H)$ generates $G\otimes H$ (see \cite[\S 59]{fuchsii}).
\epf
\mkp

For simplicity's sake, we will denote by $(G\otimes H,\otimes )$ the tensor product of two groups $G, H$ in $\sA$.

\bdfn\label{dfn001}
If $\sB\subseteq \sA$ has tensor product, we say %$(G\otimes_{\sB}H,\otimes_{\sB})$
that the \emph{tensor map} $\otimes_{\sB}$ is an \emph{epi-bihomomorphism} if for every $G,H\in\sB$ %the tensor product $(G\otimes_{\sB}H,\otimes_{\sB})$ satisfies
the following property holds: given any two group homomorphisms $h_1$ and  $h_2$, both defined on $G\otimes_{\sB}H$
and with the same range group, if
$h_1\circ \otimes_{\sB} = h_2\circ \otimes_{\sB}$ then it holds that $h_1=h_2$. It is readily seen
that, since $\otimes(G\times H)$ generates the group $G\otimes H$, it follows that $\otimes$ is an epi-bihomomorphism.
\edfn
\mkp

We now recall the definition of a reflective subcategory. Here we restrict ourselves to subcategories of $\sA$.

\bdfn\label{def 2}
Let $\mathcal{B}$ be a subcategory of $\sA$. We say that
%\begin{enumerate}
%\item
$\mathcal{B}$ is a \emph{reflective} subcategory of $\mathcal{A}$ if for each $G\in \mathcal{A}$ there is a pair $(r_\mathcal{B}G,r_\mathcal{B})$,
whith $r_\mathcal{B}G\in \mathcal{B}$ and $r_\mathcal{B}:G\rightarrow r_\mathcal{B}G$ is a continuous homomorphism, such that if $h:G\rightarrow B$ is a continuous homomorphism
with $B\in \mathcal{B}$ there exists a unique continuous homomorphism $h^r:r_\mathcal{B}G\rightarrow B$ such that the following diagram
\vspace{-0.5cm}
\begin{equation*}
\begin{picture}(148,120)
\put(40,5){$B$}
\put(-40,60){$G$}
\put(112,60){$r_\mathcal{B}G$}
\put(-15,60){\vector(1,0){121}}
\put(40,65){$r_\mathcal{B}$}
\put(-21,56){\vector(4,-3){47}}
\put(-15,30){$h$}
\put(108,56){\vector(-4,-3){47}}
\put(100,30){$h^r$}
\end{picture}
\end{equation*}
commutes.
\edfn
\mkp

When the map $r_\mathcal{B}$ is an epimorphism for every $G\in\sA$, we say that  $\mathcal{B}$ is an \emph{epireflective} subcategory of $\mathcal{A}$.
\mkp

We show below that reflective subcategories of $\sA$ have tensor product.

\bprp\label{316}  Every reflective subcategory $\sB$ of $\sA$ has tensor product $\otimes_{\sB}$. %$(G\otimes_{\sB}H,\otimes_{\sB})$.
Furthermore, if $\sB$ is epireflective,
then the bihomomorphism $$\otimes_{\sB}\colon G\times H\to G\otimes_{\sB}H$$ is an epi-bihomomorphism
and the homomorphism $b^{\otimes_{\sB}}$, that makes the diagram
\begin{equation*}
\begin{picture}(148,120)
\put(40,5){$B$}
\put(-60,60){$G\times H$}
\put(112,60){$G\otimes_{\sB}H$}
\put(-15,60){\vector(1,0){121}}
\put(40,65){$\otimes_{\sB}$}
\put(-21,56){\vector(4,-3){47}}
\put(-15,30){$b$}
\put(108,56){\vector(-4,-3){47}}
\put(100,30){$b^{\otimes_{\sB}}$}
\end{picture}
\end{equation*}

commutative, is unique.
\eprp
\pf\
Let $\sB$ be a reflective subcategory of $\sA$ and let $G,H$ be two groups in $\sB$. Given a continuous bihomomorphism
$b:G\times H\rightarrow B$ { with $B\in \mathcal{B}$, by Theorem \ref{th 1}}, there is a unique continuous homomorphism $b^\otimes\colon G\otimes H\to B$ such that
$b= b^\otimes\circ\otimes$. Now, since $\sB$ is reflective, there is a unique continuous homorphism $b^{\otimes_{\sB}}: r_\mathcal{B}[G\otimes H]\rightarrow B$ such that
$b^\otimes=b^{\otimes_{\sB}}\circ r_\mathcal{B}$. Then taking $G\otimes_{\sB}H=r_\mathcal{B}[G\otimes H]$ and  ${\otimes_{\sB}}=r_\mathcal{B}\circ\otimes$,
we obtain that $(G\otimes_{\sB}H,\otimes_{\sB})$ is a $\mathcal{B}$-tensor product of $G$ and $H$.

Now, assume that $\sB$ is epireflective and suppose $h_1$ and $h_2$ are two group homomorphisms defined on $G\otimes_{\sB}H=r_\mathcal{B}[G\otimes H]$ such that
$h_1\circ r_\mathcal{B}\circ\otimes = h_2\circ r_\mathcal{B}\circ\otimes$.
Since $\otimes(G\times H)$ generates the group $G\times H$, it follows that $h_1\circ r_\mathcal{B} = h_2\circ r_\mathcal{B}$ and,
since we are assuming that $r_\mathcal{B}$ is an epimorphisms in $\sA$, it follows that  $h_1=h_2$. This verifies that
$\otimes_{\sB}=r_\mathcal{B}\circ\otimes$ is an epi-bihomomorphism. The unicity of $b^{\otimes_{\sB}}$ is a straightforward consequence of this fact.
\epf
\mkp

\brem\label{rem1}
We notice that the above result implies that the tensor products with epi-bihomomorphisms are essentially unique.
\erem

Indeed, suppose that we have two groups $G$ and $H$ with two tensor products $(T_1,\phi_1)$ and $(T_2,\phi_2)$ for some subcategory $\sB$ of $\sA$
such that $\phi_1$ and $\phi_2$ are both epi-bihomomorphisms.
Applying the commutativity of diagram above to $\phi_2$, we have
\begin{equation*}
\begin{picture}(148,120)
\put(40,5){$T_2$}
\put(-55,60){$G\times H$}
\put(112,60){$T_1$}
\put(-15,60){\vector(1,0){121}}
\put(40,65){$\phi_1$}
\put(-21,56){\vector(4,-3){47}}
\put(-15,30){$\phi_2$}
\put(108,56){\vector(-4,-3){47}}
\put(100,30){$\widetilde{\phi_2}$}
\end{picture}
\end{equation*}
\mkp

\noindent which yields $\phi_2=\widetilde{\phi_2}\circ\phi_1$ and, in like manner, we have $\phi_1=\widetilde{\phi_1}\circ\phi_2$. Therefore
$$\phi_2=(\widetilde{\phi_2}\circ\widetilde{\phi_1})\circ\phi_2.$$
Since we also have the diagram
\begin{equation*}
\begin{picture}(148,120)
\put(40,5){$T_2$}
\put(-55,60){$G\times H$}
\put(112,60){$T_2$}
\put(-15,60){\vector(1,0){121}}
\put(40,65){$\phi_2$}
\put(-21,56){\vector(4,-3){47}}
\put(-15,30){$\phi_2$}
\put(108,56){\vector(-4,-3){47}}
\put(100,30){$\hbox{id}_{T_2}$}
\end{picture}
\end{equation*}

\noindent it follows that $\widetilde{\phi_2}\circ\widetilde{\phi_1}=\hbox{id}_{T_2}$. Analogously, we have
$\widetilde{\phi_1}\circ\widetilde{\phi_2}=\hbox{id}_{T_1}$, which implies that the groups $T_1$ and $T_2$ are topologically isomorphic.

Thus, we can unambiguously denote by $(G\otimes_{\sB}H,\otimes_{\sB})$ the tensor product of two groups $G$ and $H$
if the map $\otimes_{\sB}$ is an epi-bihomomorphism.
%\erem
\mkp

Let $\sA_2$ denote the category of Hausdorff topological abelian groups. It is readily seen that $\sA_2$ is an epireflective subcategory of $\sA$.
As a consequence, we have:

\bcor\label{cor001}
The category $\sA_2$ of all Hausdorff topological Abelian groups has a unique tensor product.
\ecor
\mkp

Next, we introduce two important classes of topological groups { (cf. \cite[Page 422, exercises 7.1.g and 7.1.h]{Arh_Tka:book_TG})}.

\bdfn\label{def 02}
A non-empty class $\sB$ of topological groups is an \emph{$\sS\sC$-variety} (resp., \emph{$\overline{\sS\sC}$-variety})
if it is closed under taking arbitrary subgroups and topological products (resp., closed under taking closed subgroups and topological products).
\edfn
%It is easily seen that {$\sS\sC$-varieties} and {$\overline{\sS\sC}$-varieties} are both reflective. Therefore, the following result holds true.

\bcor\label{Cor01}
Every subcategory $\sB$ of $\sA$ that is an $\sS\sC$-variety
(resp., every subcategory $\sB_2$ of $\sA_2$ that is an $\overline{\sS\sC}$-variety) has a unique tensor product.
\ecor
\pf\ It suffices to notice that every subcategory $\sB$ of $\sA$ that is an $\sS\sC$-variety is epireflective in $\sA$ and
 every subcategory $\sB$ of $\sA_2$ that is an $\overline{\sS\sC}$-variety is epireflective in $\sA_2$.\epf

We can now associate to each subcategory $\sB\subseteq\sA$ (resp., $\sB\subseteq\sA_2$) a category with tensor product in a canonical way.

\bprp\label{Pr 1}  { For every subcategory $\sB\subseteq\sA$ (resp., $\sB\subseteq\sA_2$) there exists a subcategory
$\sB^\otimes\supseteq \sB$ that has tensor product.}
\eprp
\pf\ It suffices to define $\sB^\otimes$ as the smallest $\sS\sC$-variety containing $\sB$ (resp.,   the smallest $\overline{\sS\sC}$-variety containing $\sB$).\epf
\mkp

Here, we are interested in an epireflective  subcategory of $\sA$ containg the sub\-ca\-te\-gory $\sL$ of locally compact abelian groups in order to
set an appropriate frame for the tensor product of two groups in $\sL$. In this sense, one might think of the subcategory of pro-Lie groups
as a good candidate for it. However, it will be shown later on that the notion of tensor product appears quite naturally as a \emph{dual group}
and { there are pro-Lie groups with bad reflexivity properties} (see \cite{FerHerSha:pams}).
This forces us to search for an epireflective subcategory that is also closed under taking dual groups. A safe option for this purpose
is the subcategory $\mathcal Q$ of Hausdorff locally quasi-convex groups that is dealt with in the next section.

\section{Locally quasi-convex groups}

 As we have pointed out { in} the Introduction, Garling \cite{Garling} { introduced groups with duality} in order to define the tensor product of two  locally compact abelian groups.
In this section, we connect Garling's definition with the more general approach we have taken here.
Thus, we provide an alternative and somewhat more general proof that the tensor product defined by Garling is topologically isomorphic to
the tensor product in the subcategory $\sQ$ of Hausdorff locally quasi-convex groups.

First, we recall some basic definitions about the duality theory of topological abelian groups
that complement the notions presented in the Introduction. All groups are assumed to be
Abelian, Hausdorff, topological groups from here on. Therefore, we will use the additive notation and $0$ for the trivial element.
In particular, we identify the unit circle group { of} the complex plane $\TT$ with the additive group $(-1/2,1/2]$,
where the addition is defined modulo $\ZZ$.
%\emph{A character} on a topological group $G$ is a continuous group homomorphism from $G$ to
%the torus group $\TT$ and the set of all characters on $G$, with pointwise addition, is a group denoted by $\widehat{G}$.

%Let \IndexPrint{$\K(G)$} denote the family of all compact subsets of $G$.
For a set $A\subseteq G$ and a positive real $\epsilon$, define
$$[A,\epsilon] = \{\chi \in \widehat{G}: (\forall a \in A)\ |\chi(a)|\leq\epsilon\}.$$\
The sets $[K,\epsilon]\subseteq\widehat{G}$, with $K$ compact subset of $G$ and $\epsilon>0$,
form a neighborhood base at the trivial character,
defining the compact-open topology.
%We also denote by {$\widehat{G}$} the \emph{topological} group obtained in this manner.

%$G$ is \emph{reflexive} if the evaluation map\
%$$\frak e :G\to\widehat{\widehat{G\,}}$$
%defined by $\mathcal E(g)(\chi)=\chi(g)$ for all $g\in G,\chi\in\widehat{G}$,
%is a topological isomorphism.
%The \emph{Pontryagin-van Kampen Theorem} asserts that every locally compact abelian group is reflexive.

Let $K$ be a compact subset of $G$. For each $n$, the set $K_n=K\cup 2K\cup\dots\cup nK$ is
compact, and $[K_n,1/4]\subseteq [K,1/4n]$. Thus, the sets $[K,1/4]$, with $K$ a compact subset of $G$, also form a neighborhood
base of $\widehat G$ at the trivial character.

\bdfn
For $A\subseteq G$, $A^\rhd = [A,1/4]$\index{$A^\rhd$}.
Similarly, for $X\subseteq\widehat{G}$, $X^\lhd =\{g\in G:(\forall \chi \in X)\ |\chi(g)|\leq \frac{1}{4}\}$.
\edfn

\blem[{\cite[Proposition 1.5]{Banaszczyk}}]\label{Upol}
For each neighborhood $U$ of $0$ in $G$, $U^\rhd$ is a compact subset of $\widehat{G}$.
\elem

\bdfn[Vilenkin \cite{Vilenkin}]
A set $A\subseteq G$ is \emph{quasi-convex} if $A^{\rhd\lhd}=A$.
$G$ is \emph{locally quasi-convex} if it has a neighborhood base at its identity consisting of quasi-convex sets.
We denote by $\mathcal Q$ the subcategory of Abelian Hausdorff locally quasi-convex groups.
\edfn
For each $A\subseteq G$, $A^\rhd$ is a quasi-convex subset of $\widehat{G}$.
Thus, $\widehat{G}$ is locally quasi-convex for all topological groups $G$.
Moreover, local quasi-convexity is hereditary for arbitrary subgroups.

$A^{\rhd\lhd}$ is the smallest quasi-convex subset of $G$ containing $A$, and is closed.

In the case where $G$ is a topological vector space, $G$ is locally quasi-convex %in the present sense
if, and only if, $G$ is a locally convex topological vector space in the ordinary sense { \cite[Proposition 2.4]{Banaszczyk}}.

If $G$ is a Hausdorff locally quasi-convex group, its characters separate points of $G$,
and thus the evaluation map $\mathcal E:G\to \widehat{\widehat G\,}$ is injective.
For each quasi-convex neighborhood $U$ of $0$ in $G$, $U^\rhd$ is a compact subset of $\widehat{G}$
(Lemma \ref{Upol}), and thus $U^{\rhd\rhd}$ is a neighborhood of $0$ in $\widehat{\widehat G\,}$.
As $\mathcal E[G]\cap U^{\rhd\rhd}=\mathcal E[U^{\rhd\lhd}]=\mathcal E[U]$, we have that $\mathcal E$ is open in $\mathcal E(G)$
(see \cite[Lemma 14.3]{Banaszczyk}).
The following fact is folklore. \medskip

\bprp\label{Prop_00001}
The subcategory {$\mathcal Q$} is epireflective.
Therefore, it has a unique tensor product.
\eprp
\medskip

We now discuss the approach of tensor product given by Garling (loc. cit.).

Given two locally quasi-convex groups $G$ and $H$, we denote by $\XX(G, H)$ the group of all continuous bicharacters of $G\times H$.
By the definition of tensor product, it is easily seen that the dual group of $G\otimes H$ (or, equivalently $G\otimes_\mathcal Q H$) is
algebraically isomorphic to the group $\XX(G,H)$. We include this fact below for future reference in the manuscript.

\bprp\label{Prop_0001}
Given two arbitrary locally quasi-convex groups $G$ and $H$, we have $$\widehat{G\otimes H}\simeq \widehat{G\otimes_\mathcal Q H}\simeq \XX(G,H).$$
\eprp

Now, endow $\XX(G, H)$ with the compact open topology; that is,
a neighborhood base at $0$ is given by the sets $[K,\epsilon]$  with $K$ a compact subset in $G\times H$ and $\epsilon>0$.
The symbol $\XX'(G,H)$ denotes the group of all continuous characters (the dual group) of $\XX(G, H)$. Following Garling \cite{Garling},
we equip the group $\XX'(G,H)$ with the topology of the uniform convergence on the equicontinuous
subsets of $\XX(G,H)$; that is, a neighborhood base at $0$ in $\XX'(G,H)$ is given by the sets  $[E,\epsilon]$ with $E$ is an { equicontinuous} set in $\XX(G,H)$.
It is easily verified that $\XX'(G,H)\in \mathcal{Q}$. %If $E$ is an equicontinouos set in $X'(G,H)$ and $\epsilon>0$,

%If $A$ is a subset of an abelian topological group $L$, let us put $A^{\triangleright}=\{\chi \in L':  |\chi(l)|\leq \frac{1}{4}\mbox{ for every } l\in L\}$. If $A\subseteq L'$, let us put $A^{\triangleleft}=\{l\in L: |\chi(l)|\leq \frac{1}{4}\mbox{ for every } \chi \in A\}$. A subset $A$ of $L$ is quasi-convex if $A^{\triangleright\triangleleft}=A$. We SAY (ya se ha usado) that $L$ is locally quasi-convex if $L$ has a local base at the neutral element consisting of quasi-convex sets.

Consider the evaluation map $\mathcal E: G\times H\longrightarrow \XX'(G,H)$, defined by $\mathcal E(x,y)(b):=b(x,y)$, for all $(x,y),\in G\times H$ and $b\in \XX(G,H)$, which is a continuous bihomomorphism. Applying the definition of tensor product, there is a unique continuous homomorphism
$$\mathcal E^\otimes: G\otimes H\longrightarrow \XX'(G,H)$$ such that the diagram

\begin{equation*}
\begin{picture}(148,120)
\put(15,3){$\XX'(G,H)$}
\put(-70,60){$G\times H$}
\put(112,60){$G\otimes H$}
\put(-15,60){\vector(1,0){121}}
\put(40,65){$\otimes$}
\put(-21,56){\vector(4,-3){47}}
\put(-15,30){$\mathcal E$}
\put(108,56){\vector(-4,-3){47}}
\put(100,30){$\mathcal E^\otimes$}
\end{picture}
\end{equation*}
\noindent commutes.
Let $T$ be the subgroup $\mathcal E^\otimes(G\otimes H)$ in $\XX'(G,H)$, equipped with the topology inherited as a subgroup of $\XX'(G,H)$.
We show below that $(T,\mathcal E)$ is the $\mathcal{Q}-$tensor product of $G$ and $H$,
which proves that Garling's definition of the tensor product of two LCA groups is a realization of
$G\otimes_\mathcal Q H$, the $\sQ$-tensor product of $G$ and $H$ (cf. \cite[Def. 3]{Garling}).

\bthm\label{ThQ} Let $G$ and $H$ be two groups in $\sQ$ and let $T$ be the subgroup $\mathcal E^\otimes(G\otimes H)$ in $\XX'(G,H)$,
equipped with the topology inherited as a subgroup of $\XX'(G,H)$. Then, $(T, \mathcal E)$ is the $\mathcal{Q}$-tensor product of $G$ and $H$.
\ethm

\pf\
Let $(A(G\times H),\sigma)$ be the abelian free topological group generated by $G\times H$. Then there is a continuous homomorphism
$$\mathcal E^\sigma: A(G\times H)\longrightarrow T\leq \XX'(G,H)$$
such that the diagram  %$\mathcal E=\mathcal E^\sigma\circ\sigma$.

\begin{equation*}
\begin{picture}(148,148)
\put(20,113){$A(G\times H)$}
\put(40,5){$T$}
\put(-60,60){$G\times H$}
\put(112,60){$G\otimes H$}
\put(-18,60){\vector(1,0){121}}
\put(20,65){$\otimes$}
\put(50,75){\color{blue}{$\mathcal E^\sigma$}}
\put(-21,68){\vector(4,3){47}}
\put(-13,90){$\sigma$}
\put(60,102){\vector(4,-3){47}}
\put(94,90){$\pi$}
\put(-21,56){\vector(4,-3){47}}
\put(-13,30){$\mathcal E$}
\put(108,56){\vector(-4,-3){47}}
\put(95,30){$\mathcal E^\otimes$}
\put(43,103){\color{blue}{\vector(0,-1){82}}}
\end{picture}
\end{equation*}
\noindent commutes.

Let $L\in \mathcal{Q}$ and $\varphi\colon G\times H \longrightarrow L$ be a continuous bihomomorphism. Then there is a continuous homomorphism $\varphi^\sigma: A(G\times H)\longrightarrow L$ such that $\varphi=\varphi^\sigma\circ\sigma$.
We claim that if $\mathcal E^\sigma(\sum\limits_{i\in F}n_iw_i)=0$ for some $\sum\limits_{i\in F}n_iw_i\in A(G\times H)$
(here $w_i\in G\times H$, $n_i\in\ZZ$ and $F$ is a finite index subset),
then $\varphi^\sigma(\sum\limits_{i\in F}n_iw_i)=0$. \mkp

Indeed, if $\mathcal E^\sigma(\sum\limits_{i\in F}n_iw_i)=0\in \XX'(G,H)$, then for each $b\in \XX(G,H)$ we have
$$0=\mathcal E^\sigma(\sum\limits_{i\in F}n_iw_i)(b)=(\sum\limits_{i\in F}n_i\mathcal E(w_i))(b)=\sum\limits_{i\in F}n_ib(w_i).$$
Moreover, since $\chi\circ\varphi\in \XX(G,H)$ for every character $\chi\in \widehat{L}$, we deduce
$$0=\sum\limits_{i\in F}n_i(\chi\circ\varphi)(w_i)=\chi(\sum\limits_{i\in F}n_i\varphi(w_i)).$$
Taking into account that $L$ is a Hausdorff locally quasi-convex group, it follows
$$0=\sum\limits_{i\in F}n_i\varphi(w_i)=\varphi^\sigma(\sum\limits_{i\in F}n_iw_i).$$
So that $\ker \mathcal E^\sigma\subseteq \ker\varphi^\sigma$.
Therefore, there is a homomorphism $$h:T\longrightarrow L$$ such that $\varphi^\sigma=h\circ\mathcal E^\sigma$.
This yields $$h(\sum\limits_{i\in F}n_i\mathcal E(w_i))=\sum\limits_{i\in F}n_i\varphi(w_i)$$
for each $$\sum\limits_{i\in F}n_i\mathcal E(w_i)=\mathcal E^\sigma(\sum\limits_{i\in F}n_iw_i)=\mathcal E^\otimes(\sum\limits_{i\in F} n_iw_i)\in T.$$

To see that $h$ is continuous it is enough to prove that $h$ is continuous at $\mathcal E(0,0)=0$. Indeed, let $V$ be a quasi-convex neighborhood of $0\in L$.
 Then there is an equicontinuous set $E_L\in \widehat{L}$ such that $V=E_L^\triangleleft$. Therefore $E:=\{\chi\circ\varphi:\chi\in E_L\}\subseteq \XX(G,H)$ is equicontinuous
 and $E^{\triangleright}$ is a neighborhood of $0$ in $X'(G,H)$. Let us see that $$h(T\cap E^\triangleright)\subseteq V.$$
We must verify that for every
$\sum\limits_{i\in F}n_i\mathcal E(w_i)\in T\cap E^\triangleright$ and $\chi\in E_L$ it holds
$$|\chi(h(\sum\limits_{i\in F}n_i\mathcal E(w_i)))|\leq \frac{1}{4}.$$

We have

\begin{equation*}
\begin{array}{l}
|\chi(h(\sum\limits_{i\in F}n_i\mathcal E(w_i)))|=|\chi(h(\mathcal E^\sigma(\sum\limits_{i\in F}(n_iw_i)))|=|\chi((h\circ\mathcal E^\sigma)(\sum\limits_{i\in F}n_iw_i))|=\\
|\chi(\varphi^\sigma(\sum\limits_{i\in F}n_iw_i))| = |\chi(\sum\limits_{i\in F}n_i\varphi(w_i))|=|\sum\limits_{i\in F}n_i(\chi\circ\varphi)(w_i)| =
|(\sum\limits_{i\in F}n_i\mathcal E(w_i))(\chi\circ\varphi)|\leq \frac{1}{4}.
\end{array}
\end{equation*}

%Therefore $$h(\sum\limits_{i\in F}\mathcal E(w_i))\in E_L^\triangleleft=V.$$
\noindent This completes the proof.
\epf

\section{Main result}
\subsection{Basic facts}
For an abelian group $G$ (resp., topological abelian group $G$), the symbol $G^*:=Hom(G,\ZZ )$ (resp., $G':=CHom(G,\ZZ )$) denotes the group of all homomorphisms
(resp., continuous homomorphisms) of $G$ into $\ZZ$ and is called the \emph{$\ZZ$-dual group} of $G$ (resp., \emph{ topological $\ZZ$-dual group} of $G$).
A group $G$ is called \emph{$\ZZ$-reflexive} if the natural map $G\rightarrow G^{**}$ is an isomorphism. For any index set $I$,
the symbol $\ZZ^{(I)}$ denotes the direct sum of copies of $\ZZ$ and $\ZZ^I$ denotes { their direct product}.
The proof of the following result can be found in \cite[Theorem 94.4]{fuchsii} or \cite[Cor. 1]{Eda:ja1983}.

\bprp
The group $Hom(\ZZ^I,\ZZ )$ is free for any index set $I$. Furthermore, if $I$ is of non measurable cardinality, then both $\ZZ^I$ and $\ZZ^{(I)}$ are $\ZZ$-reflexive
and $Hom(\ZZ^I,\ZZ)\simeq \ZZ^{(I)}$.
\eprp
\mkp

Therefore, in the remainder of the paper, we look at { the} group $\ZZ^{(I)}$ as the $\ZZ$-dual of the group $\ZZ^{I}$,
assuming that $I$ has non measurable cardinality from here on.
Equipped with the topology, $t_p(\ZZ^{I})$, of pointwise convergence on $\ZZ^{I}$,
the group $\ZZ^{(I)}$ becomes a topological group, which is closed as a subgroup of $\ZZ^{\ZZ^I}$.
The main goal of this section is to find its Pontryagin dual.
%First, we shall make use of these groups in order to give an example of a topological group that has no predual group.

%\bprp
%The group  $\ZZ^{^(\NN)}$, equipped with the pointwise convergence $t_p(\ZZ^{\NN})$ has no predual group.
%\eprp
%\pf\ Assume there is a topological group $G$ such that $\widehat{G}\cong \ZZ^{(\NN)}$. Equip $G$ with the topology of
%pointwise convergence on the elements of $\ZZ^{^(\NN)}$, that is the topology $t_p(\ZZ^{^(\NN)})$.
%\epf

The following result was proved by Garling (loc. cit.)

\bthm\label{Garling Th5}
If $G_1$ and $G_2$ are Hausdorff locally quasi-convex groups, and if $\widehat{G_2}$ is equipped  with the topology of uniform convergence on the compact sets of $G_2$,
then there is a canonical \emph{1} to \emph{1} homomorphism of $\X(G_1, G_2)$ into $CHom(G_1, \widehat G_2)$.
If $CHom(G_1, \widehat G_2)$ is given the topology of uniform convergence
on the compact sets of $G_1$, this mapping is continuous and { open onto its image}. Further, if
$G_2$ is locally compact, the mapping is onto, so that $\X(G_1, G_2)\cong CHom(G_1,\widehat{G_2})$, algebraically and topologically.
\ethm
\mkp

\brem\label{Rem001}
The isomorphism dealt with in Theorem \ref{Garling Th5} is naturally defined as follows:
$\Phi\colon \X(G_1, G_2)\to CHom(G_1, \widehat G_2) :\ \Phi(\beta)=\tilde{\beta}$ where
$\tilde{\beta}(g_1)[g_2]:=\beta(g_1,g_2)$.
\erem

\subsection{The dual group of $\ZZ^{(I)}$}

According to Proposition \ref{Prop_0001} and Theorem \ref{Garling Th5}, we have that $\X(\ZZ^I, \TT)$,  the dual group of the tensor product $\ZZ^I\otimes_\sQ \TT$, is topologically isomorphic to
$H=CHom(\ZZ^I,\ZZ)$. However, while Garling equips the group $H$ with the compact open topology,
remark that $H$ is equipped with the topology of pointwise convergence on $\ZZ^I$ in this paper.
Taking this fact into account, we first prove that $\ZZ^I\otimes_\sQ \TT$ is canonically embedded in $\widehat{H}$. %From here on, if $\cg\in H$ and $\cx\in \ZZ^I$, we set
%$$\langle\cg,\cx\rangle :=\sum\limits_{n\in\NN}\cg(n)\cdot\cx(n).$$

\bprp \label{Prp1}
$\ZZ^I\otimes_\sQ \TT\subseteq \widehat{H}$.
\eprp
\pf\
Let $\theta\colon \ZZ^I\otimes_\mathcal Q \TT\to \widehat H$ be the map defined as follows:
$$\theta((\cx\otimes t))[{\cg}]:= t\,\cdot\, <\cg,\cx >,\ \overline x\otimes t\in \ZZ^I\otimes_\mathcal Q\TT,\
\cg \in H,$$

where  $$\langle\cg,\cx\rangle :=\sum\limits_{i\in I}\cg(i)\cdot\cx(i).$$

Since $H$ is equipped with the topology of pointwise convergence on $\ZZ^I$, it follows that the map $\cg\mapsto\ <\cg,\cx>$ of $H$ into $\ZZ$
is continuous for all $\cx\in \ZZ^I$.  On the other hand $\theta(\cx\otimes_\sQ t)$ is the composition of $<\cg,\cx>$ with the map
$n\mapsto t\cdot n$, which is obviously continuous as $\ZZ$ is discrete. This yields the continuity of $\theta(\cx\otimes_\sQ t)$.
Furthermore, it is easily seen that $\theta$ defines a $1$-to-$1$ group homomorphism and, therefore, embeds canonically
$\ZZ^I\bigotimes_\sQ \TT$ in $\widehat{H}$.
\epf

\brem\label{Rem002}
Notice that ${\widehat{H}}$ is a subgroup of $\widehat{H_d}=\TT^I$, where $H_d$ denotes the same algebraic group $H=\ZZ^{(I)}$ equipped with the discrete topology.
Thus, with some notational abuse, we may claim that $$\ZZ^I\otimes_\sQ \TT\leq \widehat{H}\leq\TT^I.$$
The question is which elements $\ct\in\TT^I$ belong to $\widehat{H}$. { In order to answer this question, we need a number of preliminary lemmas.}
\erem
\mkp

%For simplicty's sake, from here on, in this section, we denote the group $\ZZ^I$ by $X$.
\mkp

\blem\label{Lem0}
The group $H$ has a neighborhood basis of the identity { consisting of all finite intersections} $$V:=\bigcap\limits_{1\leq j\leq m} V_j,$$ where $V_j:=\{\cg\in H :\ <\cg,\cx_j>=0$\}, $\cx_j\in \ZZ^I$, $1\leq j\leq m.$
\elem
\pf\ It is enough to observe that the group $H$ is equipped with the topology that it inherits from $\ZZ^{\ZZ^I}$ and $\ZZ$ is a discrete group.
\epf

\brem\label{Rem0}
Plainly, any of these subsets $V$ is a nontrivial open subgroup, which means that $H$ is a non discrete, non archimedean, subgroup
(see \cite[Prop. 2.3]{FerHerSha:pams}, where the assertion is proved for $I=\NN$).
\erem

\blem\label{Lem1}
Let $\ct$ be an element of $\TT^I$. Then $\ct$ belongs to $\widehat{H}$ if, and only if, there is $\{\cx_1,\dots ,\cx_m\}\subseteq \ZZ^I$ \st\
if $\cg\in H$ and $<\cg,\cx_j>=0$, for all $1\leq j\leq m$, then $<\cg,\ct>=0$.
\elem
\pf\ \emph{Sufficiency} is clear.

\noindent \emph{Necessity:} Suppose that  $\ct\in \TT^I\cap \widehat{H}$ and set $U:=\TT\cap R$, where $R=[-1/4,1/4]$. %is the right half plane of the complex field.
Clearly $U$ is a neighborhood of the identity in $\TT$ that contains no nontrivial subgroups. Since we are assuming that $\ct$ is continuous, there must be a basic open subgroup $$V:=\bigcap\limits_{1\leq j\leq m} V_j,$$ {where} $$V_j:=\{\cg\in H :\ <~\cg,\cx_j>~=~0~\},\ \cx_j\in \ZZ^I,\ 1\leq j\leq m,$$
such that $\ct[V]\subseteq U$. Since $\ct[V]$ is a subgroup of $\TT$, it follows that $\ct[V]=\{0\}$. This completes the proof.
\epf
\mkp

Therefore, { for each element $\ct\in \TT^I\cap \widehat{H}$}, there is a finite subset $F:=\{\cx_1,\dots ,\cx_m\}\subseteq \ZZ^I$ satisfying the assertion formulated
in Lemma \ref{Lem1}.
There is no loss of generality in assuming that $F$ is minimal under inclusion. We say that $F$ is a  \emph{continuity subset} for $\ct$.

\blem\label{Lem2}
If\, $\ct\in \widehat{H}$ and $\ct(i)\not=0$ then, for each continuity subset $F$ for $\ct$, there is $\cx\in F$ \st\ $\cx(i)\not=0$.
\elem
\pf\ Suppose that $\cx(i)= 0$ for all $\cx\in F$. If we take $\cg\in H$ \st\ $\cg(i)=1$ and $\cg(j)=0$ if $j\not=i$, then  $<\cg,\cx>=0$ for all $\cx\in F$
and $<\cg, \ct>=\ct(i)\not=0$. This is a contradiction that completes the proof.\epf

In the sequel, if $J$ is a subset of $I$, we define $\delta_J\in\ZZ^I$ as the characteristic function of $J$; that is
$\delta_J(i)=1$ if $i\in J$ and $\delta_J(i)=0$ if $i\notin J$. Furthermore, if $\overline{x}$ and $\overline{t}$ belong to
$\ZZ^I$ and $\TT^I$, respectively, then $\delta_J\cdot\overline{x}$ and, respectively, $\delta_J\cdot\overline{t}$ are defined
as $(\delta_J\cdot\overline{x})(i)=\cx(i)$, resp., $(\delta_J\cdot\overline{t})(i)=\ct(i)$, if $i\in J$, and
$(\delta_J\cdot\overline{x})(i)=0$, resp., $(\delta_J\cdot\overline{t})(i)=0$, if $i\notin J$.
Remark that $$<\cg,\delta_J\cdot\overline{x}>=<\delta_J\cdot\cg,\overline{x}>\ \forall \cg\in H, \cx\in \ZZ^I.$$

\blem\label{Lem3}
Let $J$ be an infinite subset of $I$ and set $H_J=\{ \cg\in H : \cg(i)=0,\ \forall i\notin J\}$. Then for every open subgroup $V$ in $H$,
which is defined by a finite subset $F=\{\cx_1,\dots ,\cx_m\}\subseteq \ZZ^I$,
it holds $H_J\cap V\not=\{0\}$.
\elem
\pf\ Let $V$ be a canonical open neighborhood of $0$ in $H$. We have
$V:=\bigcap\limits_{1\leq j\leq k} V_j$, $V_j:=\{\cg\in H :\ <\cg,\cx_j>=0\}$, $\cx_j\in F$, $1\leq j\leq k$.
Set $W_j:=\{\cg\in H :\ <\cg,\delta_J\cdot \cx_j>=0\}$, $\cx_j\in F$, $1\leq j\leq k$ and
$W:=\bigcap\limits_{1\leq j\leq k} W_j$. We have %For every $\cg\in W$, we have
$$0=<\cg,\delta_J\cdot\overline{x}>=<\delta_J\cdot\cg,\overline{x}>\ \forall \cg\in W, \cx\in F,$$
which means that $\cg\in W$ if and only if $\delta_J\cdot \cg\in V$.
On the other hand, observe that there is a canonical topological isomorphism of the group $H_J$ onto $H^{(J)}$
when the latter group is equipped with the pointwise convergence on $\ZZ^J$ and, by Remark \ref{Rem0},
we know that the group $H^{(J)}$, equipped with the pointwise convergence on $\ZZ^J$, is not discrete.
Since $W\cap H_J$ is a neighborhood of $0$ in $H_J$, it follows that
$W\cap H_J\not=\{0\}$. Finally, observe that $W\cap H_J\subseteq V$, which yields $V\cap H_J\not=\{0\}$.\epf

\blem\label{Lem3.1}
The map $p_J\colon H\to H$ defined by $p_J(\cg):=\delta_J\cdot\cg$ is continuous and open { onto its image} for each subset $J$ of $I$.
\elem
\pf\ The proof is clear if $J$ is finite. Therefore, we can assume that $J$ is infinite without loss of generality.
We have seen in Lemma \ref{Lem3} that for every canonical open neighborhood $V$ of $0$, there is another one $W$ \st\
$\cg\in W$ if and only if $\delta_J\cdot \cg\in V$, which means that $W = p_J^{-1}(V\cap H_J)$. This completes the proof.\epf

\blem\label{Lem3.2}
Let $J$ be a subset of $I$ and let $\overline{t}\in\TT^I$ be an element of $\widehat{H}$.
Then the point  $\overline{s}=\delta_J\cdot\overline{t}$, defined by
$\overline{s}(i)=\overline{t}(i)$ if $i\in J$ and $\overline{s}(i)=0$ if $i\notin J$, also belongs to $\widehat{H}$.
\elem
\pf\ It suffices to observe that $\overline s=\ct\circ p_J$.

\blem\label{Lem4}
If $\ct\in \TT^I$ contains infinitely many  irrationals that are linearly independents with respect to the field of
rational numbers $\QQ$, then $\ct\notin \widehat{H}$.
\elem
\pf\ Let $J$ be an infinite subset of $I$ such that $\ct(i)\notin \QQ$ for all $i\in J$ and the set $\{\ct(i) : i\in J\}$ is { linearly independent over} $\QQ$.
Reasoning by contradiction, suppose that $\ct\in\widehat{H}$. By Lemma \ref{Lem1}, there is $\{\cx_1,\dots ,\cx_m\}\subseteq \ZZ^I$ \st\
if $<\cg,\cx_j>=0$, $1\leq j\leq m$, then $<\cg,\ct>=\{0\}$. Let $V$ be the open subgroup in $H$ defined by
$\{\cx_1,\dots ,\cx_m\}$ and take $0\not=\cg\in V\cap H_J$. Then $\sum\limits_{i\in J} \cg(i)\ct(i)=0$,
which yields $\cg(i)=0$ for all $i\in J$. This is a contradiction that completes the proof.\epf

\blem\label{Lem5}
Let $\ct\in \TT^I$ \st\ $\ct(i)=\frac{p_i}{q_i}t_0$, where $\frac{p_i}{q_i}\in\QQ$, the pair of numbers
$p_i$ and $q_i$ are { relatively prime} for all $i\in I$, and $t_0$ is an irrational or $1$.
Then $\ct\in\widehat{H}$ \sii\ the sequence $(|q_i|)$ is bounded.
\elem
\pf\ \emph{Sufficiency:} Assume that $(|q_i|)$ is bounded and, therefore, defines a finite set. Let $M:=lcm (|q_i|)$,
where the acronym lcm stands for \emph{least common multiple}, and set $\overline{x}:=(\frac{p_iM}{q_i})\in \ZZ^I$.
If $$\langle \cg,\overline{x} \rangle=0,$$ then
$$0=\langle \cg,(\frac{p_iM}{q_i}) \rangle=M\langle \cg,(\frac{p_i}{q_i}) \rangle,$$
which implies $$\langle \cg,(\frac{p_i}{q_i}) \rangle=0.$$ This means that
$$\langle \cg,(\frac{p_i}{q_i}t_0) \rangle =0$$ and as a consequence $(\frac{p_i}{q_i}t_0)\in\widehat{H}$.

\emph{Necessity:}
assume that $(\frac{p_i}{q_i}t_0)\in\widehat{H}$. By Lemma \ref{Lem1}, there is an open neighborhood of the
identity $V$ defined by a continuity set $\{\cx_1,\dots ,\cx_m\}\subseteq \ZZ^I$ \st\
if $\langle\cg,\cx_j\rangle=0$, $1\leq j\leq m$, then $\langle\cg,(\frac{p_i}{q_i}t_0)\rangle=\{0\}$.
\mkp

\noindent CLAIM 1\label{claim1} We may assume that the set
$\{\cx_1,\dots ,\cx_m\}$ is linearly independent on the field $\QQ$ without loss of generality.

\noindent\emph{Proof of Claim 1.} Suppose that $\cx_1=\gl_2\cx_2+\dots +\gl_m\cx_m$,
where $\gl_j=a_j/b_j$, $\{a_j,b_j\}\subseteq \ZZ$, $2\leq j\leq m$.
Set $M:=lcm(\{|b_2|,\dots ,|b_m|\})$. Then $M\cx_1=M\gl_2\cx_2+\dots +M\gl_m\cx_m$ and
$M\gl_j\in\ZZ$ for $2\leq j\leq m$. Therefore, if $\langle\cg,\cx_j\rangle=0$, $2\leq j\leq m$,
then $\langle\cg,M\cx_1\rangle=\{0\}$. Hence $\langle\cg,\cx_1)\rangle=(1/M)\langle\cg,M\cx_1)\rangle=\{0\};$
i.e., we can dispose of $\cx_1$ without modifying the definition of $V$.
This completes the proof of Claim~1.\mkp

Define the map $$\Phi\colon\ZZ^{(I)}\to \ZZ^m$$ by
$$\Phi(\cg):=(\langle\cg,\cx_1\rangle,\dots , \langle\cg,\cx_m\rangle)$$
and extend it to $$\Phi\colon\QQ^{(I)}\to \QQ^m$$ with the same definition:
$$\Phi((x_i)):=(\langle(x_i),\cx_1\rangle,\dots , \langle(x_i),\cx_m\rangle).$$

Set $$\phi\colon\QQ^{(I)}\to \QQ$$ by  $$\phi((x_i)):=\langle(x_i),(\frac{p_i}{q_i})\rangle\in\QQ.$$
\mkp

\noindent CLAIM 2\label{claim2}
$\ker\, \Phi\subseteq \ker\, \phi$.

\noindent\emph{Proof of Claim 2.} Take $(y_i)\in \ker\,\Phi$. Since $(y_i)\in \QQ^{(I)}$, there is $m\in\NN$ \st\
$m(y_i)\in \ZZ^{(I)}$. Furthermore $\Phi(m(y_i))=m\Phi((y_i))=0$. Hence
$\langle m(y_i),(\frac{p_i}{q_i}t_0)\rangle=0$, which implies $\phi(m(y_i))=0$. Therefore
$\phi((y_i))=(1/m)\phi(m(y_i))=0$. This completes the proof of Claim 2.
\mkp

\noindent CLAIM 3\label{claim3}
There is a linear map $\overline{\phi}\colon\QQ^m\to\QQ$ \st\ $\overline{\phi}\circ\Phi=\phi$.

\noindent\emph{Proof of Claim 3.}
Since the set $\{\cx_1,\dots ,\cx_m\}$ is linearly independent on the field $\QQ$ and minimal as a continuity set of $V$, it follows that
for every $\overline{z}\in\QQ^m$, there is $(y_i)\in\QQ^{(I)}$ \st\ $\Phi((y_i))=\overline{z}$.
Define $\overline{\phi}(\overline{z}):=\phi((y_i))$. By Claim \ref{claim2}, the map $\overline{\phi}$ is well defined
and it is readily seen that satisfies the assertion of Claim 3.\mkp

Now, we have $$\overline{\phi}(\overline{z})=\sum_{j=1}^{m}\gl_jz_j,$$
$$\Phi((y_i))=(\langle(y_i),\cx_1\rangle,\dots , \langle(y_i),\cx_m\rangle)=
(\sum\limits_i y_i\cx_1(i),\dots , \sum\limits_i y_i\cx_m(i))$$ and
$$\phi((y_i))=\langle(y_i),(\frac{p_i}{q_i})\rangle =\sum\limits_i y_i\cdot (p_i/q_i).$$

Since $\phi((y_i))=\overline{\phi}(\Phi((y_i))$, we conclude
$$\frac{p_i}{q_i}=\sum_{j=1}^{m} \gl_j\cx_j(i),\ \forall i\in\NN.$$

Here, $\cx_j(i)\in\ZZ$ and $\gl_j=a_j/b_j\in\QQ$. Set $M:= lcm(\{|b_1|,\dots , |b_m|\}$. Then
$$\frac{p_iM}{q_i}=\sum_{j=1}^{m} M\gl_j\cx_j(i)\in\ZZ,\ \forall i\in\NN.$$

Since $p_i$ and $q_i$ are { relatively prime}, this means that $M$ is a multiple of $|q_i|$ for all $i\in\NN$,
which means that the sequence $(|q_i|)$ is bounded.\epf
\mkp

We can now prove the main result in this section.

\bthm\label{theorem1}
$\widehat{H}\cong \ZZ^I\otimes_\sQ \TT$.
\ethm
\pf\ Let $\ct\in \TT^I$ be an arbitrary element of $\widehat{H}$. According to Lemma \ref{Lem4},
the set $\{\ct(i) : i\in I\}$ { can contain at most } finitely many irrationals that are linearly independent over $\QQ$.
Thus, there are:
\begin{enumerate}
  \item a finite set $\{t_0,\dots , t_m\}$, where $t_0$ can possibly be $1$ and the subset $\{t_1,\dots , t_m\}$
consists of irrational numbers that are linearly independent over $\QQ$,
  \item a finite partition $\{I_0,\dots, I_m\}$ of $\NN$, and
  \item a sequence $(\frac{p_i}{q_i})$, where
{ the
numbers} $p_i$ and $q_i$ are { relatively prime} for all $i\in\NN$,
\end{enumerate}

\noindent \st\ for every $1\leq j\leq m$ and $i\in I_j$, it holds $\ct(i)=\frac{p_i}{q_i}t_j$.
Since $\ct\in\widehat{H}$, by Lemma \ref{Lem3.2}, we know that $\ct_j:=\delta_{I_j}\cdot \ct$ also
belongs to $\widehat{H}$, $0\leq j\leq m$. Now, according to Lemma \ref{Lem5}, the sequence
$(|q_i|)_{i\in I_j}$ is bounded and, as a consequence, finite. That is, each $q_i$ with $i\in I_j$
can only take a finite number of values $\{q_{1,j},\dots, q_{m_j,j}\}$. Set $I_{ij}:=\{k\in I_j: q_k=q_{ij}\}$.
Then the point $\ct_{ij}\in \TT^I$, defined by $\ct_{ij}(k)=p_k\cdot\frac{t_j}{q_{ij}}$ if $k\in I_{ij}$
and $\ct_{ij}(k)=0$ if $k\notin I_{ij}$, belongs to $\ZZ^I\otimes_\sQ \TT$. Hence
$\ct_j=\sum\limits_{1\leq i\leq m_j}\ct_{ij}$ is also an element of $\ZZ^I\otimes_\sQ \TT$ for all $0\leq j\leq m$.
This means that $t=\sum\limits_{0\leq j\leq m}\ct_{j}$ belongs to $\ZZ^I\otimes_\sQ \TT$, which completes the proof.\epf
\mkp

We are now in position to prove our main result.

\textbf{{ Proof of Theorem} \ref{Thm_A}.}
It suffices to take $I=\NN$ in Theorem \ref{theorem1}.
\epf

\section{The dual group of $A_p(X)$}

The aim of this section is to find the dual group of $A_p(X)$, that is the free abelian group $A(X)$ equipped with
the topology $t_p(C(X,\ZZ))$. First, we need to recall some topological notions. A \emph{$0$-dimensional} space is a topological space having an open { basis}
consisting of \emph{clopen} subsets. The weight ${w}(X)$ of a topological space $X$ is the minimum cardinality of an open basis of $X$.
A { topological} space is said to be $\NN$-compact if it can be embedded as a closed subset in a product of copies of $\NN$. Clearly,
every $\NN$-compact space is $0$-dimensional. A function $f$ defined on a topological space $X$ is called \emph{$k$-continuous} if
$f_{|K}$ is continuous { on} $K$ for each compact subspace $K$ of $X$. For { a} $0$-dimensional space $X$,
define a space $k_\NN X$ to be the space { underlying} $X$ equipped with the smallest topology making each $\ZZ$-valued, $k$-continuous
function on $X$ continuous. The space $X$ is said to be a \emph{$k_\NN$-space} if $k_\NN X=X$.

%It is a well known fact that for any topological space $X$, there is a $0$-dimensional
%space $\tau_\NN X$ for which $C(X,\ZZ)\simeq C(\tau_\NN X,\ZZ)$. Furthermore,
{ A topological space is called \emph{$\NN$-compact} if it is homeomorphic to a closed subspace of the product of copies of the countable discrete space $\NN$.
It is a well known fact that for any $0$-dimensional space $X$}, there exists a unique
$\NN$-compactification $\gb_\NN X$ of $X$ to which each element of $C(X,\ZZ)$ admits a continuous extension, and then
$C(X,\ZZ)\simeq C(\gb_\NN X,\ZZ)$. In view of these facts, we may { restrict} out attention to $0$-dimensional spaces or even
$\NN$-compact spaces.
%Along this section $X$ designates a $0$-dimensional $\NN$-compact space.

Our first result is a direct consequence of Proposition \ref{Prop_0001} and Theorem \ref{Garling Th5}. %recall some results about groups of $\ZZ$-valued continuous fucntions.

\bprp \label{Prp41}
Let $X$ be a $0$-dimensional space. { Then}
$$(C_p(X,\ZZ)\otimes_\mathcal{Q}\TT)\ha\simeq\XX(C_p(X,\ZZ),\TT)\simeq CHom(C_p(X,\ZZ),\ZZ).$$
\eprp

\blem\label{Lem42}
Let $X$ be a $0$-dimensional space. Then $$CHom(C_p(X,\ZZ),\ZZ)\simeq A(X).$$
\elem
\pf\ Let $\phi\in CHom(C_p(X,\ZZ),\ZZ)$. We will see that there is a minimum finite subset $F\subseteq X$ such that
for $f\in C(X,\ZZ)$,
if $f(F)=\{0\}$ then $\phi(f)=0$. We call this subset $F$ the \emph{support} of $\phi$, written ${\rm supp}(\phi)$.

Indeed, assuming that is not true, for every finite subset $F\subseteq X$, there is $f_F\in C_p(X,\ZZ)$ such that
$f_F(F)=\{0\}$ and $\phi(f_F)\not=0$. %In fact, there is no lost of generality in assuming $\phi(f)=1$.
It is easily seen that the net $\{f_F : F\ \text{is a finite subset of}\ X\}$, ordered by inclusion, converges to $0$ in $C_p(X,\ZZ)$.
However $\phi(f_F)\not=0$ for all finite subset $F\subseteq X$, which is a contradiction. { Thus, there is some finite subset $F\subseteq X$
\st\ if $f(F)=\{0\}$ then $\phi(f)=0$. Let $\mathcal F$ be the collection of all finite subsets in $X$ with this property. We claim that
if $D,E$ belong to $\mathcal{F}$ then so does $D\cap E$. Indeed, suppose that $f\in C(X,\ZZ)$ and $f(D\cap E)=\{0\}$. If $E\subseteq D$ then $D\cap E=E$, which belongs to $\mathcal F$.
Therefore $\phi(f)=0$. Hence, we can safely  assume that $E\setminus D\not=\emptyset$. Let $U$ be a clopen subset of $X$ such that
$E\setminus D\subseteq U$ and $U\cap D=\emptyset$. Clearly $f=f\chi_U+(1-\chi_U)f$, where $\chi_U$ denotes the characteristic function
of $U$. Now, since $U\cap D=\emptyset$,  we have  $(f\chi_U)(D)=0$
and, since $f(E\cap D)=\{0\}$, we have $f(1-\chi_U)(E)=0$. Hence $\phi(f)=\phi(f\chi_U)+\phi(f(1-\chi_U))=0$.
Therefore, since the family $\mathcal F$ consists of
finite subsets, it follows that there is a minimum subset $F$ satisfying the property
formulated above.}

Let $F$ be the support of $\phi$. Then there is a homomorphism $\tilde{\phi}\colon\ZZ^F\to \ZZ$ such that
the following diagram

\vspace{-0.5cm}
\begin{equation*}
\begin{picture}(148,120)
\put(40,5){$\ZZ$}
\put(-70,60){$C_p(X,\ZZ)$}
\put(112,60){$\ZZ^F$}
\put(-15,60){\vector(1,0){121}}
\put(40,65){$r_F$}
\put(-21,56){\vector(4,-3){47}}
\put(-15,30){$\phi$}
\put(108,56){\vector(-4,-3){47}}
\put(100,30){$\tilde{\phi}$}
\end{picture}
\end{equation*}
commutes.

From this fact, it follows that $\phi(f)=\sum\limits_{x\in F} n_xf(x)=\sum\limits_{x\in F} n_x\delta_x(f)$.
The map $$\Delta\colon CHom(C_p(X,\ZZ),\ZZ)\to A(X)$$ defined by $$\Delta(\phi)= \sum\limits_{x\,\in\, {\rm supp}(\phi)} n_x\delta_x$$
is clearly an algebraic onto isomorphism. This completes the proof.\epf

\bcor\label{Cor43}
Let $X$ be a $0$-dimensional space. { Then}
$$(C_p(X,\ZZ)\otimes_\mathcal{Q}\TT)\ha\simeq A(X).$$
\ecor
\mkp

In the same way as in Proposition \ref{Prp1}, we obtain

\blem\label{Lema44}
$C_p(X,\ZZ)\otimes_\sQ \TT\subseteq \widehat{A_p(X)}$.
\elem
\pf\
Let the map $\theta\colon C_p(X,\ZZ)\otimes_\sQ \TT\to \widehat{A_p(X)}$ defined as follows:
let $f\otimes t$ be an arbitrary element of $C_p(X,\ZZ)\otimes_\sQ\TT$ where $f \in C_p(X,\ZZ)$ and $t\in\TT$ and
let $\phi$ be an element of $A_p(X)$.  Then
$$\theta((f\otimes_\sQ t))[{\phi}]:= t\,\cdot\, \phi(f),$$
where  $$\phi(f) :=\sum\limits_{x\,\in\, {\rm supp}(\phi)} n_x\delta_x(f)=\sum\limits_{x\,\in\, {\rm supp}(\phi)} n_xf(x).$$
Since $A_p(X)$ is equipped with the topology of pointwise convergence on $C_p(X,\ZZ)$, it follows that the map $\phi\mapsto\ \phi(f)$ of
$A_p(X)$ into $\ZZ$ is continuous for all $f\in C_p(X,\ZZ)$.  On the other hand $\theta(f\otimes_\sQ t)$ is the composition of $\phi(f)$
with the map $n\mapsto t\cdot n$, which is obviously continuous as $\ZZ$ is discrete. This yields the continuity of $\theta(f\otimes_\sQ t)$.
Furthermore, it is easily seen that $\theta$ defines a $1$-to-$1$, group isomorphism, therefore, embeds canonically
$C_p(X,\ZZ)\bigotimes_\sQ \TT$ into $\widehat{A_p(X)}$.
\epf
\mkp

\brem\label{Rem45}
The group $C_p(X,\ZZ)$ is a subgroup { of} $\ZZ^X$ and therefore $t_p(C_p(X,\ZZ)\subseteq t_p(\ZZ^X)$. Hence
$\widehat{A_p(X)}\subseteq \ZZ^X\otimes_\sQ \TT$. %The follwing lemma clarifies which elements of $\ZZ^X$
%define continuous characters of $A_p(X)$ when $X$ is a $k_\NN$-space.
\erem
\mkp

We are now in position to prove the main result in this section.
\mkp

\textbf{{ Proof of Theorem} \ref{Thm_B}.}

It will suffice to prove that for every $0$-dimensional, $k_\NN$-space $X$, if $f\in\ZZ^X$ satisfies that $f\otimes_\mathcal{Q}t\in \widehat{A_p(X)}$,
then there exists $g\in C(X,\ZZ)$ such that $f\otimes_\mathcal{Q}t=g\otimes_\mathcal{Q}t$.

%\blem\label{Lem46}
%Let $X$ be a $0$-dimensional, $k_\NN$-space and let $f\otimes_\mathcal{Q}t\in \widehat{A_p(X)}$. Then $g\in C(X,\ZZ)$ such that
%$f\otimes_\mathcal{Q}t=g\otimes_\mathcal{Q}t$.
%\elem
%\pf\
First, notice that, since $X$ is $0$-dimensional, it follows that $t_p(C(X,\ZZ))_{|X}$ induces the original topology on $X$.

Let $K\subseteq X$ be a compact set. Since $f\otimes_\mathcal{Q}t$ is continuous on $K$, it follows that
$f\otimes_\mathcal{Q}t(K)=\{t\cdot f(x) : x\in K\}$ is compact for all $t\in\TT$.
This means that the subset $f(K)\subseteq \ZZ$ is compact in the Bohr topology of $\ZZ$. Now, it is a well known fact
that every Bohr-compact subset in $\ZZ$ is finite (see \cite{GalHdez:Fundamenta}), which means that
$f(K)$ must be finite for every compact subset $K$ of $X$.

Now, assume first that $t$ is an irrational number in $\TT$ and
set $K_n:=\{x\in K : f(x)=n\}$. We shall prove that the collection $\{K_n : n\in f(K)\}$ forms a finite partition of $K$
by clopen subsets.

Indeed, if $x\in\overline{K_n}$, then $t\cdot f(x)=t\cdot n$ in $\TT$, which implies that
$t\cdot f(x)- t\cdot n\in\ZZ$ and, since $t$ is an irrational in $\TT$, it follows that $f(x)=n$.
Therefore, each subset $K_n$ is closed and therefore clopen in $K$.

Thus, we have obtained that the collection $\{f^{-1}(n) : n\in f(K)\}$ forms a finite partition of $K$
by clopen subsets. This means that $f$ is $k$-continuous and since $X$ is a $k_\NN$-space, it follows
that $f\in C(X,\ZZ)$. In this case, we can take $g=f$.

In case $t$ is a rational number in $\TT$, take $m\in\NN$ the smallest natural number such that $t\cdot m\in\ZZ$ and
set $X_n:=\{ x\in X : f(x)=km+n\ \hbox{for some}\ k\in\ZZ\}$ for each natural number $n$ with $0\leq n<m$.
Notice that if $x\in X_n$, we have
$f\otimes_\mathcal{Q}t(x)=t\cdot (km+n)=t\cdot n$ in $\TT$. We claim that every subset $X_n$ is closed.
Indeed, if $x\in\overline{X_n}$, then $t\cdot (f(x)-n)=0$ in $\TT$. Therefore, there is
$k\in\ZZ$ such that $f(x)-n=km$, which implies that $x\in X_n$.
Therefore, we have obtained that the collection $\{X_n : 0\leq n< m\}$ forms a finite partition of $X$
by clopen subsets. Set $g\colon X\to\ZZ$ defined by $g_{|X_n}=n$. Then $g$ is clearly
continuous and $f\otimes_\mathcal{Q}t=g\otimes_\mathcal{Q}t$. \epf

\section{Groups of { integer-valued} continuous functions}

In this section, we apply our previous results in order to study the reflexivity properties of the abelian groups
$C_p(X,\ZZ)$ consisting of all continuous functions from $X$ into the discrete additive group of the integers $\ZZ$.
%\bdfn\label{Def5.1} (Def5.1)
%A topological space $X$ is \emph{$0$-dimensional} if it contains an open basis consisting of clopen subsets.
%The space $X$ is called \emph{$\NN$-compact} if it is homemorphic to  closed subspace of a product $\NN^\kappa$
%for some cardinal $\kappa$. A $\ZZ$-valued function on a space $X$ is said to be \emph{$k$-continuous} if the restricition
%$f_{|K}$ is continuous for each compact subset of $K$. Given a topological space $X$, we denote by $k_\ZZ X$
%the set $X$ equipped with the smallest topology so that every $\ZZ$-valued $k$-continuous function on $X$ is continuous.
%It is known that when $X$ is an $\NN$-compact space then $C_p(X,\ZZ)'\simeq C_p(X,\ZZ)^*$ (see \cite[Cor. 6.7]{Edo_Kiyosawa_Ohta:1989}).
%\edfn
%\mkp
From here on, all spaces are assumed to be $0$-dimensional, $\NN$-compact and of non-measurable cardinality.
The first result in this section is straightforward but we formulate it below for future reference.

\bprp\label{Prp51}
Let $X$ be a $0$-dimensional { space}. Then $$C_p(X,\ZZ)\ha\simeq (\ZZ^X)\ha\simeq\TT^{(X)}\simeq A(X)\otimes_{\mathcal{Q}} \TT.$$
\eprp
\pf\ The first isomorphism follows from the fact that $C_p(X,\ZZ)$ is dense in $\ZZ^X$. The other isomorphisms are folklore.
\epf
\mkp

In what follows we investigate when the group $C_p(X,\ZZ)$ is Pontryagin reflexive. For that purpose, we refer to \cite{Her_Usp:JMAA2000}
where the reflexivity of the group $C_p(X,\RR)$ is investigated in depth. %There, the following result is proved.
We next prove that the group $C_p(X,\ZZ)$ is dually embedded in $C_p(X,\RR)$.
From here on, we use multiplicative notation in order to better explain the key ideas { of the proofs}.

\bprp\label{Prp52}
Let $\chi\colon C_p(X,\ZZ)\to \TT$ be a continuous character. Then there is $\tilde{\chi}\colon C_p(X,\RR)\to \TT$
such that $\tilde{\chi}_{|C_p(X,\ZZ)}=\chi$.
\eprp
\pf\ %Here, we use multiplicative notation in order to better explain the key ideas of the proof.
Since $\chi\in \TT^{(X)}$, we have that $\chi=\sum\limits_{x\in X} t_x\delta_{x}$, where $t_x\in\TT$ and the sum is finite.
Therefore, for every $f\in C(X,\ZZ)$ we have $\chi(f)=\prod\limits_{x\in X} t_x^{f(x)}$.
For each $t_x\not=1$, there is a real number $\lambda_x\in (0,1)$ such that $\exp(2\pi i\lambda_x)=t_x$. Then we have
$$\chi(f)=\prod \exp(2\pi i\lambda_x)^{f(x)}=\prod \exp(2\pi i\lambda_x f(x))=\exp(2\pi i\sum \lambda_x f(x)).$$
Define $\tilde{\chi}(f):= \exp(2\pi i\sum \lambda_x f(x))$ for all $f\in C(X,\RR)$. It is clear that $\tilde{\chi}$
is a continuous character on $C_p(X,\RR)$ that extends $\chi$. \epf
\mkp

The dual space of $C_p(X)$ can be identified with the space $L(X)$ of
all formal linear combinations of points of $X$ with real coefficients.
If $a=\sum \lambda_x x\in L(X)$, the finite set $\{x\in X:\lambda_x\ne0\}$ is called
the {\it support\/} of $a$ and is denoted by $\supp\, a$.
If $B\subseteq L(X)$, we set $\supp\, B=\bigcup \{\supp\, a:a\in B\}$.
\mkp

\brem\label{Rem53}
We notice that for { a  $0$-dimensional space} $X$, the following isomorphisms hold:
$$C_p(X,\RR)\ha\simeq (L(X),t_c(C_p(X,\RR))\simeq (\RR^{(X)}, t_c(C_p(X,\RR))$$ and
$$C_p(X,\ZZ)\ha\simeq (A(X)\otimes_{\mathcal{Q}} \TT,t_c(C_p(X,\ZZ))\simeq (\TT^{(X)}, t_c(C_p(X,\ZZ)).$$

The first chain of isomorphism is well known (see \cite{Her_Usp:JMAA2000}) and the second one is consequence
of Proposition \ref{Prp51}.
\erem
\mkp

Next result is a variation of \cite[Proposition 3.1]{Her_Usp:JMAA2000}. Here, the symbolism $(\TT^{(X)}, t_c(C(X,\ZZ))$
denotes the space $\TT^{(X)}$ equipped with the compact open topology on $C(X,\ZZ)$.

\blem\label{Lem53}
Let $X$ be a $0$-dimensional space. Then every compact subset in  $(\TT^{(X)}, t_c(C(X,\ZZ))$ has finite support.
\elem
\pf\ %Again, we use multiplicative notation in order to better explain the key ideas of the proof.
Let $D$ be a subset of $\TT^{(X)}$ whose support is
infinite. It will suffice to prove that $D$ is not countably compact.
In the sequel, we denote each element $w\in \TT^{(X)}$ as a finite sum $w=\sum_{x\in X} t_{x}$ with $t_{x}\in \TT$.

According to \cite[Lemma 3.2]{Her_Usp:JMAA2000},
there exists an infinite collection $\sV$ of pairwise disjoint
open sets in $X$ each of which meets $\supp\, D$. Construct by induction
a sequence $\{w_n=\sum_{x_n\in X} t_{x_n}\}_{n<\omega }\subseteq D$
and a sequence $\{V_n\}_{n<\omega }\subseteq \sV$
such that
every $V_n$ meets $\supp\, w_n$ and does not meet $\supp\, w_j$ whenever
$j<n$. We claim that the sequence  $\{w_n\}_{n<\omega }$ has no accumulation point in $\TT^{(X)}$.

Indeed, let $w_0$ be an arbitrary but fixed element in $\TT^{(X)}$. %and
%assume without loss of generality that $(\supp\, w_0) \cap V_n =\emptyset$ for all $n<\omega$.
For every $n<\omega$, pick $x_n\in (\supp\, w_n)\cap V_n$ and a clopen
neighbourhood $W_n$ of $x_n$ such that $W_n\subseteq V_n$ and
$(\supp\, w_n)\cap W_n=\{x_n\}$. %and $(\supp\, w_0)\cap W_n=\emptyset$.

For any $n<\omega$, pick $g_n\in C(X,\ZZ)$ such that $\Re((t_{x_n})^{g_n(x_n)})<0$
(the real part of $(t_{x_n})^{g_n(x_n)}$ is less than $0$)
and $g_n(X\setminus W_n)=\{0\}$. The sequence
$\{g_n\}_{n<\omega }$ converges pointwise to zero,
since the cozero-set of each $g_n$ is contained in $W_n\subseteq V_n$
and the sequence $\{V_n\}_{n<\omega }$ is disjoint.

%Let $w_0$ be an arbitrary element of $\TT^{(X)}$.
Since the family $\sV$ is disjoint,
we have that $\supp w_0$ can meet finitely many member of $\sV$ at most.
So there is $N\in\NN$ such that $(\supp w_0)\cap W_n=\emptyset$ for every $n>N$.
The set
$K_N=\{g_n\}_{N<n<\omega }\cup\{0\}$ is compact in $C_p(X,\ZZ)$.
Since $g_n(w_0)=1$ for all $n>N$, it follows that $w_0\in (K_N)^\rhd$.
However, since $\Re((t_{x_n})^{g_n(x_n)})<0$ for all $n\in\NN$, it follows that
$\{w_n\}_{n<\omega }\cap (K_N)^\rhd\subseteq \{w_1,\dots w_N\}$, a finite set.
Thus $w_0$ may not be an accumulation point of $\{w_n\}_{n<\omega }$.
This completes the proof. \epf
\mkp

\bprp\label{Prp54}
The evaluation map $\mathcal E :C_p(X,\ZZ)\rightarrow C_p(X,\ZZ)\haha$ is a
topological isomorphism of $C_p(X,\ZZ)$ onto its image.
\eprp
\pf\
Since the group $C_p(X,\ZZ)\leq \ZZ^X$, it follows that $C_p(X,\ZZ)$ is locally quasi-convex. Therefore,
the evaluation is open in its image. On the other hand, in order to verify the continuity of $\mathcal E$,
it suffices to prove that
every compact subset $K$ of
$(\TT^{(X)}, t_c(C_p(X,\ZZ))=C_p(X,\ZZ)\ha$\, is equicontinuous.

Lemma \ref{Lem53} implies that $K$ is contained in a subgroup of the form
of the form
$K_{F}=\oplus_{x\in F}\TT_x$, where $F$ is a finite subset of $X$.
Plainly $K_{F}$ is equicontinuous, hence so is $K$.\epf

\textbf{{ Proof of Theorem} \ref{Thm_C}.}

Let $X$ be a $0$-dimensional space such that $C_p(X,\RR)$ is reflexive, which we know exists by
Theorem 3.10 and Example 3.11 in \cite{Her_Usp:JMAA2000}. By Proposition \ref{Prp54},
it follows that the evaluation mapping is continuous and open in its image. Therefore, { in} order to verify that $C_p(X,\ZZ)$ is reflexive,
it will suffice to prove that for every $0$-dimensional space $X$ it holds $$C_p(X,\ZZ)\ha\ha\simeq C_p(X,\RR)\ha\ha\cap\ZZ^X.$$

Set $$\pi\colon \RR^{(X)}\to \TT^{(X)}$$ the canonical quotient mapping defined by
$$\pi(\sum_{x\in X} \lambda_x)= \exp(2\pi i\sum_{x\in X} \lambda_x).$$ We notice that $\ker\pi=\ZZ^{(X)}$.

Take an arbitrary element $f$ in $C_p(X,\RR)\ha\ha\cap\ZZ^X$. We have that $$f\colon (\RR^{(X)}, t_c(C_p(X,\RR))\to \TT$$ is continuous and
$f_{|\ZZ^{(X)}}=0$. This means that $f$ factors through a continuous map $\tilde{f}$ defined on $\TT^{(X)}$  such that the
following diagram
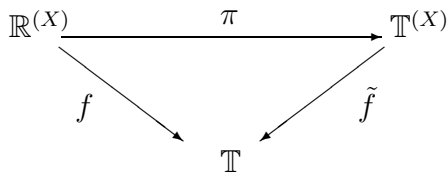
\vspace{-0.5cm}
\begin{equation*}
\begin{picture}(148,120)
\put(40,9){$\TT$}
\put(-40,60){$\RR^{(X)}$}
\put(104,60){$\TT^{(X)}$}
\put(-20,60){\vector(1,0){121}}
\put(40,65){$\pi$}
\put(-21,56){\vector(4,-3){47}}
\put(-15,30){$f$}
\put(102,56){\vector(-4,-3){47}}
\put(92,30){$\tilde{f}$}
\end{picture}
\end{equation*}
commutes.

Set $$\Phi\colon C_p(X,\RR)\ha\ha\cap\ZZ^X\to C_p(X,\ZZ)\ha\ha$$ defined by $\Phi(f)=\tilde{f}$.
{ It is easily seen} that $\Phi$ is an onto topological group isomorphism.\epf

\noindent \textbf{Acknowledgements:} { The authors wish to thank the referee for her/his constructive report that has helped us improve this article.}

\end{document}